\documentclass{article}
\usepackage[T2A]{fontenc}
\usepackage[utf8]{inputenc}
\usepackage[english,russian]{babel}
\usepackage{amssymb,amsmath,amsthm,amsfonts}

\begin{document}

\begin{center}\bf

On planar harmonic functions with identical Jacobians

\end{center}

{\bf Annotation:} This note provides a complete description of a family of sense preserving harmonic functions in the open unit disk $\mathbb{D}$ that have the same Jacobians, provided that one of the representatives of this family is known.

\bigskip


\bigskip

{\bf 1. Введение.} Напомним, что функция $f$, определеная и дважды непрерывно дифференцируемая в области $D\subset\mathbb{C}$, называется гармонической, если всюду в $D$ имеет место тождество Лапласа $f_{z\overline{z}}(z)\equiv0$. Нижними индексами $z$ и $\overline{z}$ здесь и далее обозначается соответственно применение к $\mathbb{R}$-дифференцируемой комплекснозначной функции следующих дифференциальных операторов 
\[ \frac{\partial}{\partial z}=\frac{1}{2}\left(\frac{\partial}{\partial x}-i\,\frac{\partial}{\partial y}\right), \hspace{1cm} \frac{\partial}{\partial \overline{z}}=\frac{1}{2}\left(\frac{\partial}{\partial x}+i\,\frac{\partial}{\partial y}\right). \]
Хорошо известно, что всякая гармоническая в $D$ функция $f$ может быть представлена в виде суммы $f=h+\overline{g}$, где $h$, $g$ -- голоморфные в $D$ функции. Функцию $h$ называют голоморфной частью функции $f$, а $\overline{g}$ -- антиголоморфной частью $f$. Якобиан $J$ функции $f$, которая может рассматриваться, как отображение из $\mathbb{R}^2$ в $\mathbb{R}^2$, имеет вид $J=|h'|^2-|g'|^2$. Всюду ниже речь будет идти о сохраняющих ориентацию гармонических в области функциях. По теореме Леви, для того, чтобы гармоническая функция $f$ сохраняла ориентацию в $D$, необходимо и достаточно, чтобы при любом $z\in \mathbb{D}$ было выполнено неравенство $J(z)>0$. Для таких функций определяется еще одна важная характеристика -- аналитическая дилатация $\omega=g'/h'$. Используя аналитическую дилатацию, можно сформулировать условие сохранения ориентации, эквивалентное приведенному выше: $|\omega|<1$ в $D$. Теорему Леви, а также другие классические результаты теории однолистных гармонических функций на комплексной плоскости могут быть найдены в монографии \cite{dur1}.

Напомним также \cite{dur2, hille}, что производная Шварца локально однолистной голоморфной функции $f$, определенной в $D$, определяется формулой
$$
S[f, z]=\left(\frac{f''}{f'}\right)'-\frac{1}{2}\left(\frac{f''}{f'}\right)^2.
$$
Известно, что производные Шварца двух локально однолистных голоморфных функций $f$ и $g$ совпадают тогда и только тогда, когда $g=L\circ f$, где $L$ -- невырожденное дробно-линейное преобразование плоскости $\mathbb{C}$.

Данная заметка является естественным развитием работы \cite{gn}. Главным результатом в ней является теорема, описывающая свойства якобианов сохраняющих оиентацию гармонических в $\mathbb{D}={|z|<1}$ функций, а также содержащая достаточные условия, при выполнении которых положительная функция $J\in C^{\infty}(\mathbb{D})$ является якобианом некоторой сохраняющей оиентацию гармонической в $\mathbb{D}$ функции.

\medskip
{\bf Теорема А.}
\textsl{ Положительная функция $J\in C^{\infty}(\mathbb{D})$ является якобианом некоторой сохраняющей ориентацию гармонической в $\mathbb{D}$ функции $f=h+\overline{g}$ тогда и только тогда, когда существует такая голоморфная функция $\omega$, $|\omega|<1$, что}

\begin{equation} \label{ln_main}
(-\ln\,J)_{z\overline{z}}=|\omega'|^2(1-|\omega|^2)^{-2}.
\end{equation}

\textsl{ Более того, данная функция $J$ и соответствующая ей функция $\omega$ обладают следующими свойствами: }

1. \textsl{либо функция $\ln\,J$ определена и является гармонической в $\mathbb{D}$, а функция $\omega$ постоянна;}

2. \textsl{либо функция $R=\ln(-J^2\,(\ln\,J)_{z\overline{z}})$ определена и является гармонической в  $\mathbb{D}\setminus Z$, где $Z$ -- множество изолированных нулей функции $(-\ln\,J)_{z\overline{z}}$. В этом случае $Q=2\left(   \ln(-\ln\,J)_{z\overline{z}} \right )_{zz}-\left(   \ln(-\ln\,J)_{z\overline{z}}  \right)_{z}^2$ голоморфна в $\mathbb{D}\setminus Z$, а функция $\omega$ является решением дифференциального уравнения $2S[\omega, z]=Q,$ определенного в $\mathbb{D}\setminus Z$.}

\medskip

Якобианы сохраняющих ориентацию гармонических в $\mathbb{D}$ функций, удовлетворяющие условию 1 теоремы А, далее будем называть \textsl{якобианами первого типа}, а остальные -- \textsl{якобианами второго типа}.

Кроме того, в \cite{gn} была описана структура множества сохраняющих ориентацию гармонических в $\mathbb{D}$ функций с одинаковыми якобианами второго типа в терминах их аналитических дилатаций.

\medskip
{\bf Теорема Б.}
\textsl{ Пусть $f_0=h_0+\overline{g_0}$ -- сохраняющая ориентацию гармоническая в $\mathbb{D}$ функция с аналитической дилатацией $\omega_0=g'_0/h'_0$ и якобианом $J_0$ второго типа. Тогда гармоническая в $\mathbb{D}$ функция $f$ имеет якобиан $J_0$ тогда и только тогда, когда ее аналитическая дилатация имеет вид $\omega=T\circ\omega_0$, где $T$ -- дробно-линейный автоморфизм круга $\mathbb{D}$.}

\medskip

Наконец, в работе \cite{gn} были полностью описаны гармонические функции с якобианами первого типа. 

 \medskip
{\bf Теорема В.} 
1) \textsl{Любая функция вида $f=h+a\cdot\overline{h}+b$, где $h$ -- локально однолистная голоморфная в $\mathbb{D}$ функция, $|a|<1$, $b\in\mathbb{C}$, является сохраняющей ориентацию гармонической в $\mathbb{D}$ функцией с якобианом первого типа. Обратно, всякая сохраняющая ориентацию гармоническая в $\mathbb{D}$ функция с якобианом первого типа может быть представлена в данном виде.}

2) \textsl{Пусть $f$ -- сохраняющая ориентацию гармоническая в $\mathbb{D}$ функция с якобианом первого типа. Тогда множество сохраняющих ориентацию гармонических в  $\mathbb{D}$ функций с тем же якобианом исчерпывается функциями вида $\hat{f}(z)=e^{i\alpha}h+e^{i\beta}a\overline{h}+b,$ где $h$ -- аналитическая часть $f$, $\alpha, \beta \in \mathbb{R}$, $a=f_{\overline{z}}(0)/f_z(0)$, $b\in\mathbb{C}$.}
\medskip

Заметим, что эти теоремы были доказаны в случае, когда областью определения $D$ гармонических функций является $\mathbb{D}$. Поскольку интерес пока вызывают односвязные $D$, достаточно рассмотреть случаи $D=\mathbb{D}$ и $D=\mathbb{C}$. В \cite{gn} доказано, что сохраняющие ориентацию гармонические функции, определенные в $\mathbb{C}$, имеют якобианы первого типа и, следовательно, попадают под условия теоремы В. Поэтому далее тоже будем считать, что $D=\mathbb{D}$.

Ниже доказан аналог теоремы В для случая гармонических функций с якобианами второго типа.

\bigskip

{\bf 2. Оcновной результат.} Рассмотрим прежде всего важный частный случай, который позволит решить поставленную задачу.

Пусть $v$ -- непостоянная голоморфная в $\mathbb{D}$ функция, $|v|<1$. Рассмотрим  $J=1-|v|^2$. 

Прямыми вычислениями проверяется, что функция $R$ (см. теорему А), соответствующая построенной функции $J$, является гармонической функцией в $\mathbb{D}\setminus Z$, где $Z$ -- множество нулей голоморфной функции $v'$. 

Рассуждениями, приведенными в \cite[теорема 1]{gn}, доказывается, что функция $Q=2\left(   \ln(-\ln\,J)_{z\overline{z}} \right )_{zz}-\left(   \ln(-\ln\,J)_{z\overline{z}}  \right)_{z}^2$ определена и голоморфна в $\mathbb{D}\setminus Z$. В данном случае $Q$ является удвоенной производной Шварца от функции $v$, т.е. $ Q=  2S[v, z].$

Итак, в $\mathbb{D}\setminus Z$ определено дифференциальное уравнение
$S[\omega, z]=S[v, z]$ В силу приведенного во введении свойства производной Шварца, получаем, что $\omega$ имеет вид $$ \omega=\frac{a\, v+b}{c\, v+d}, \hspace{1cm} ad-bc\neq0. $$ Среди всех этих $\omega$ равенству \eqref{ln_main} удовлетворяют лишь $\omega=T\circ v$, где $T$ -- дробно-линейный автоморфизм круга $\mathbb{D}$, что следует из теоремы Б.

Таким образом, по теореме А, функция $J=1-|v|^2$ порождает семейство сохраняющих ориентацию гармонических в $\mathbb{D}$ функций с якобианом $J$. Аналитические компоненты функций $f=h+\overline{g}$ данного семейства являются решениями системы
$$
 \begin{cases}
   |h'|^2-|g'|^2=1-|v|^2,\\
   \dfrac{|g'|}{|h'|}=\left|\dfrac{v+z_0}{1+\overline{z_0}\,v}\right|,
 \end{cases}
$$
где $z_0\in\mathbb{D}$ -- параметр. Решая систему, получаем семейство функций следующего вида
$$f_{\alpha,\beta,z_0,C}=\dfrac{e^{i\alpha}}{\sqrt{1-|z_0|^2}}\left(z+\overline{z_0}\int\,v\,dz\right)+\dfrac{e^{i\beta}}{\sqrt{1-|z_0|^2}}\overline{\left(z_0\,z+\int v\,dz\right)}+C, 
$$
где $\alpha,\beta\in\mathbb{R}, z_0\in\mathbb{D}, C\in\mathbb{C}$. Для простоты записи в дальнейшем нижние индексы, показывающие зависимость от параметров, мы будем опускать.

\bigskip

Перейдем теперь к рассмотрению общего случая.

Пусть $f=h+\overline{g}$ -- произвольная сохраняющая ориентацию гармоническая в $\mathbb{D}$ функция с якобианом $J=|h'|^2-|g'|^2$ второго типа и аналитической дилатацией $\omega=g'/h'$. Действуя так же, как в рассмотренном выше случае, построим семейство функций, порожденное якобианом $\widetilde{J}=1-|\omega|^2$. Очевидно, производные аналитических компонент гармонических функций данного семейства будут иметь следующий вид:
$$
\widetilde{h}'=\dfrac{e^{i\alpha}}{\sqrt{1-|z_0|^2}}\left(1+\overline{z_0}\,\omega\right), \hspace{0.25cm} \widetilde{g}'=\dfrac{e^{i\beta}}{\sqrt{1-|z_0|^2}}\left(\omega+z_0\right), \hspace{0.25cm} \alpha,\beta \in\mathbb{R}, z_0\in\mathbb{D}.
$$
При различных фиксированных значениях параметров $\alpha,\beta \in\mathbb{R}, z_0\in\mathbb{D}$ будем домножать $\widetilde{h}'$ и $\widetilde{g}'$ на $h'$ и затем интегрировать. Тем самым, меняя параметры, получим семейство $\cal{J}$ функций вида

\begin{equation}\label{main}
F=\dfrac{e^{i\alpha}}{\sqrt{1-|z_0|^2}}\left(h+\overline{z_0}\, g\right)+\dfrac{e^{i\beta}}{\sqrt{1-|z_0|^2}}\overline{\left(g+z_0\,h\right)}+C,
\end{equation}
где $\alpha,\beta\in\mathbb{R}, z_0\in\mathbb{D}, C\in\mathbb{C}.$

После элементарных преобразований \eqref{main} может быть переписано в более удобном виде:

\begin{equation}\label{main2}
F=A_{\alpha,\beta,z_0}\circ R_{\alpha,\beta}[f],
\end{equation}
где $R_{\alpha,\beta}[f]=e^{i\alpha}h+e^{i\beta}\overline{g}$, а $A_{\alpha,\beta,z_0}(z)=(1-|z_0|^2)^{-1/2}(z+e^{i(\alpha+\beta)}\overline{z_0\,z})+C$.

Теперь все готово для того, чтобы сформулировать и доказать следующую теорему.

\medskip

{\bf Теорема.} \textsl{Пусть $f$ -- сохраняющая ориентацию гармоническая в $\mathbb{D}$ функция с якобианом $J$ второго типа, $\alpha,\beta\in\mathbb{R}, z_0\in\mathbb{D}, C\in\mathbb{C}$ -- произвольные параметры. Рассмотрим семейство $\cal{J}$ функций, порожденное $f$. Тогда}

1) \textsl{$f\in\cal{J};$}

2) \textsl{Каждая функция семейства $\cal{J}$ имеет якобиан $J$;}

3) \textsl{Семейство $\cal{J}$ содержит все сохраняющие ориентацию гармонические в $\mathbb{D}$ функции с якобианом $J$.}

\smallskip

{\bf Доказательство.} Для доказательства 1) достаточно в \eqref{main2} положить $\alpha=\beta=z_0=C=0$.

2) Пусть $F\in\cal{J}$. Очевидно, что якобианы функций $R_{\alpha,\beta}[F]$ и $A_{\alpha,\beta,z_0}(z)$ соответственно тождественно равны $J$ и $1$. Но тогда $J_F\equiv J$, как якобиан композиции двух функций.

3) Пусть $F$ -- произвольная сохраняющая ориентацию гармоническая в $\mathbb{D}$ функция  с якобианом $J$ и аналитической дилатацией $\omega_F$. Покажем, что найдутся такие $\alpha,\beta\in\mathbb{R}, z_0\in\mathbb{D}, C\in\mathbb{C}$, что $F$ представима в виде \eqref{main}.

В силу теоремы Б имеем $$ \omega_F\equiv e^{i\gamma}\dfrac{\omega+z_0}{1+\overline{z_0}\,\omega} $$ при некоторых значениях параметров $\gamma\in\mathbb{R}$ и $z_0\in\mathbb{D}$.

Аналитические компоненты функции $F=H+\overline{G}$ являются решениями системы

$$
 \begin{cases}
   |H'|^2-|G'|^2=|h'|^2-|g'|^2,\\
   \dfrac{|G'|}{|H'|}=\left|\dfrac{\omega+z_0}{1+\overline{z_0}\,\omega}\right|.
 \end{cases}
$$
Решая эту систему относительно $|H'|$ и $|G'|$, получаем
\begin{equation}\label{H}
|H'|^2=\dfrac{|h'|^2-|g'|^2}{1-\left|\dfrac{\omega+z_0}{1+\overline{z_0}\,\omega}\right|^2}
=\dfrac{\left(|h'|^2-|g'|^2\right)|1+\overline{z_0}\,\omega|^2}{|1+\overline{z_0}\,\omega|^2-|\omega+z_0|^2}
=\dfrac{|h'+\overline{z_0}\,g'|^2}{1-|z_0|^2},
\end{equation}

\begin{equation}\label{G}
|G'|=\left|\dfrac{\omega+z_0}{1+\overline{z_0}\,\omega}\right|\,\dfrac{|h'+\overline{z_0}\,g'|}{\sqrt{1-|z_0|^2}}
=\dfrac{|\omega+z_0|}{|1+\overline{z_0}\,\omega|}\,\frac{|h'|\,|1+\overline{z_0}\,\omega|}{\sqrt{1-|z_0|^2}}
=\dfrac{|g'+z_0\,h'|}{\sqrt{1-|z_0|^2}}.
\end{equation}
Из \eqref{H} и \eqref{G} следует требуемое. {\it Теорема доказана.}

\medskip

{\bf Благодарности.} Автор данной работы является стипендиатом Фонда развития теоретической физики и математики «БАЗИС».

\end{document}